\documentstyle[amssymb,amsfonts]{amsart}

 at 10 true pt

\def\R{{\hbox{\bf R}}}

\def\allt#1{%
\smash{
\vtop{%
     \ialign{%
        ##\crcr
        $\hfil\displaystyle{\tilde \forall}\hfil$\crcr%
        \noalign{\kern1.5pt\nointerlineskip}
        $\hfil\!\!#1\hfil$\crcr\noalign{\kern1.5pt}
        }
       }
      } \hbox{$\vphantom{#1}$}
     }

\def\be#1{\begin{equation}\label{#1}}
\def\bas{\begin{align*}}
\def\eas{\end{align*}}
\def\bi{\begin{itemize}}
\def\ei{\end{itemize}}

\def\Z{{\hbox{\bf Z}}}

\newenvironment{proof}{\noindent {\bf Proof} }{\endprf\par}
\def \endprf{\hfill  {\vrule height6pt width6pt depth0pt}\medskip}
\def\emph#1{{\it #1}}
\def\textbf#1{{\bf #1}}

\theoremstyle{plain}
  \newtheorem{theorem}[subsection]{Theorem}
  \newtheorem{conjecture}[subsection]{Conjecture}

\theoremstyle{remark}

\theoremstyle{definition}

\include{psfig}

\begin{document}

\title[Fuglede's conjecture fails in 5 dimensions]{Fuglede's conjecture is false in 5 and higher dimensions}

\author{Terence Tao}
\address{Department of Mathematics, UCLA, Los Angeles CA 90095-1555}
\email{tao@@math.ucla.edu}


\begin{abstract}  We give an example of a set $\Omega \subset \R^{5}$ which is a finite union of unit cubes, such that $L^2(\Omega)$ admits an orthonormal basis of exponentials $\{ \frac{1}{|\Omega|^{1/2}} e^{2\pi i \xi_j \cdot x}: \xi_j \in \Lambda \}$ for some discrete set $\Lambda \subset \R^{5}$, but which does not tile $\R^{5}$ by translations.  This answers (one direction of) a conjecture of Fuglede \cite{fuglede} in the negative, at least in 5 and higher dimensions.
\end{abstract}

\maketitle

\section{Introduction}

Let $\Omega$ be a domain in $\R^n$, i.e., $\Omega$ is a Lebesgue measurable subset of $\R^n$ with finite non-zero
Lebesgue measure. We say that a set $\Lambda \subset \R^n$ is
a {\it spectrum} of $\Omega$ if
${\{\frac{1}{|\Omega|^{1/2}} e^{2\pi i x \cdot \xi} \}}_{\xi \in \Lambda}$ is an orthonormal basis of $L^2(\Omega)$.  In this paper we use $|\Omega|$ to denote the Lebesgue
measure of a set $\Omega$, and $\# A$ to denote the cardinality of a finite set $A$.

\begin{conjecture}\cite{fuglede} A domain $\Omega$ admits a spectrum if and only if it is possible to tile $\R^n$ by a family of 
translates $\{t + \Omega: t \in \Lambda \}$ of $\Omega$ (ignoring sets of measure zero). 
\end{conjecture}
 
Fuglede \cite{fuglede} proved this conjecture (also known as the spectral set conjecture) under the additional assumption that the 
tiling set or the spectrum are lattice subsets of $\R^n$.  This conjecture arose from the study of commuting self-adjoint extensions of the partial derivative operators $\frac{\partial}{\partial x_j}$, and has attracted much recent interest, see the references given in the bibliography for a partial list of papers relating to this conjecture, and \cite{Kol3} for a survey. 

Our main result here is that Fuglede's conjecture is false in sufficiently high dimension.

\begin{theorem}\label{thm:main}  Let $n \geq 5$ be an integer.  Then there exists a compact set $\Omega_2 \subset \R^n$ of positive measure such that $L^2(\Omega_2)$ admits an orthonormal basis of exponentials $\{ \frac{1}{|\Omega_2|^{1/2}} e^{2\pi i \xi_j \cdot x}: \xi_j \in \Lambda_2 \}$ for some $\Lambda_2 \subset \R^n$, but such that $\Omega_2$ does not tile $\Z^n$ by translations.  In particular, Fuglede's conjecture is false in $\R^n$ for $n \geq 5$.
\end{theorem}

The counterexample $\Omega_2$ is elementary - it is an explicit finite union of unit cubes, and is based on a counterexample to Fuglede's conjecture in a specific finite abelian group.  Basically, the idea is to exploit the existence of Hadamard matrices (i.e. orthogonal matrices whose entries are all $\pm 1$, or more generally a $p^{th}$ root of unity) of order not equal to a power of $p$; when $p=2$ the first example occurs at dimension 12, and when $p=3$ the first example occurs at dimension 6.  Such Hadamard matrices quickly lead to a counterexample to Fuglede's conjecture in the finite groups $\Z_2^{12}$ and $\Z_3^6$ (actually to $\Z_2^{11}$ and $\Z_3^5$), and one can use standard transference techniques to move this counterexample to $\Z^5$ and thence to $\R^5$.

Our arguments do not preclude the possibility that the conjecture may still be true in lower dimensions, and in particular in one dimension; see for instance \cite{laba} for some evidence in favor of the one-dimensional conjecture.  For instance, the results of \cite{laba} show that Fuglede's conjecture is true for cyclic $p$-groups $\Z_{p^N}$, so one cannot directly replicate the above counterexample in one dimension.

It may still be true that Fuglede's conjecture still holds in higher dimensions under more restrictive assumptions on the domain $\Omega_2$, for instance if one enforces convexity (our example is highly non-convex, although it does not fall into the class of non-convex objects studied on \cite{KP}, or the near-cubic objects studied in \cite{KL}). 

We do not address the issue as to whether the converse direction of Fuglede's conjecture might still hold; in other words, whether every set which tiles $\R^n$ by translations admits a spectrum.  Again, it seems that one should first look at $p$-groups to determine the truth or falsity of this conjecture.

\section{The finite model: failure of Fuglede in $\Z_2^{12}$, $\Z_3^6$, and $\Z_3^5$.}

We begin with a finite version of Theorem \ref{thm:main}, in the finite group $\Z_p^n$, where $p=2,3$ and $\Z_p := \Z/p\Z$.  If $x = (x_1, \ldots, x_n)$, $\xi = (\xi_1, \ldots, \xi_n)$ are elements of $\Z_p^n$, we define the dot product $\xi \cdot x \in \Z_p$ as
$$ \xi \cdot x := \sum_{j=1}^{12} \xi_j x_j$$
and in particular we can define the quantity $e^{2\pi i (\xi \cdot x)/p}$, which is always a $p^{th}$ root of unity.

To illustrate the method, we begin with a counterexample in $\Z_2^{12}$, although we will not directly use this example for our main result.

\begin{theorem}\label{thm:finite}  There exists a non-empty subset $\Omega_0 \subset \Z_2^{12}$ such that $l^2(\Omega_0)$ admits an orthonormal basis of exponentials $\{ \frac{1}{(\# \Omega_0)^{1/2}} e^{2\pi i (\xi_j \cdot x) / 2}: \xi_j \in \Lambda_0 \}$ for some $\Lambda_0 \subset \Z_2^{12}$, but such that $\Omega_0$ does not tile $\Z_2^{12}$ by translations.
\end{theorem}

\begin{proof}
Let $e_1, \ldots, e_{12}$ be the standard basis for $\Z_2^{12}$, thus $e_j$ is the 12-tuple which equals 1 in the $j^{th}$ entry and 0 everywhere else.  We shall take $\Omega_0$ to simply be this 12-element set:
$$ \Omega_0 := \{e_1, e_2, \ldots, e_{12} \}.$$
It is clear that $\Omega_0$ does not tile $\Z_2^{12}$ by translations, since $\# \Omega_0 = 12$ does not divide evenly into $\# \Z_2^{12} = 2^{12}$.  On the other hand, $\Omega_0$ admits an orthonormal set of exponentials.  To see this, we take any Hadamard matrix $H$ of order 12, for instance\footnote{We are grateful to Neil Sloane for making this example available on his web page {\tt http://www.research.att.com/$\sim$njas/hadamard}.}
$$ H :=
\left( \begin{array}{llllllllllll}
+1 & -1 & -1 & -1 & -1 & -1 & -1 & -1 & -1 & -1 & -1 & -1 \\
+1 & +1 & -1 & +1 & -1 & -1 & -1 & +1 & +1 & +1 & -1 & +1 \\
+1 & +1 & +1 & -1 & +1 & -1 & -1 & -1 & +1 & +1 & +1 & -1 \\
+1 & -1 & +1 & +1 & -1 & +1 & -1 & -1 & -1 & +1 & +1 & +1 \\
+1 & +1 & -1 & +1 & +1 & -1 & +1 & -1 & -1 & -1 & +1 & +1 \\ 
+1 & +1 & +1 & -1 & +1 & +1 & -1 & +1 & -1 & -1 & -1 & +1 \\
+1 & +1 & +1 & +1 & -1 & +1 & +1 & -1 & +1 & -1 & -1 & -1 \\
+1 & -1 & +1 & +1 & +1 & -1 & +1 & +1 & -1 & +1 & -1 & -1 \\
+1 & -1 & -1 & +1 & +1 & +1 & -1 & +1 & +1 & -1 & +1 & -1 \\
+1 & -1 & -1 & -1 & +1 & +1 & +1 & -1 & +1 & +1 & -1 & +1 \\
+1 & +1 & -1 & -1 & -1 & +1 & +1 & +1 & -1 & +1 & +1 & -1 \\
+1 & -1 & +1 & -1 & -1 & -1 & +1 & +1 & +1 & -1 & +1 & +1 
\end{array}\right).$$
One can verify that all the rows of $H$ are orthogonal to each other (this is part of what it means for a matrix to be Hadamard).  We then define the set of 12 frequencies $\Lambda_0 := \{\xi_1, \ldots, \xi_{12}\}$ by requiring that the 12-dimensional vector $(e^{2\pi i (e_j \cdot \xi_k)/2})_{j=1}^{12}$ matches the $k^{th}$ row of $H$, thus for instance
$$ \xi_1 := (1,0,0,0,0,0,0,0,0,0,0,0), \xi_2 := (1,1,0,1,0,0,0,1,1,1,0,1), \hbox{ etc.}$$
It is then clear that the twelve exponential functions $\frac{1}{\sqrt{12}} e^{2\pi i (\xi_k \cdot x)/2}$, $k = 1, \ldots, 12$, form an orthonormal basis of $l^2(\Omega)$ as claimed.
\end{proof}

The above simple example may be compared with the example used to disprove Tijdeman's conjecture in \cite{LS}, or to disprove Keller's conjecture in \cite{lagarias}; the three counterexamples are not directly related to each other, but they do share a similar flavor, in that the combinatorics of tiling and orthogonality in ``higher-dimensional'' situations may behave quite differently from what one might intuitively extrapolate from ``low-dimensional'' examples.  These three examples also show that the main obstructions to tiling or spectral conjectures in Euclidean spaces in fact come from finite abelian groups of relatively small order.

Now let $\omega := e^{2\pi i/3}$ be a cube root of unity.  Observe that
the $6 \times 6$ matrix
$$ H :=
\left( \begin{array}{llllll}
1   &   1   &   1    &   1    &   1   &   1   \\
1   &   1   &\omega&\omega&\omega^2 & \omega^2  \\
1   &\omega & 1 & \omega^2 & \omega^2   & \omega   \\
1   &\omega &\omega^2 & 1   & \omega  &  \omega^2   \\
1   &\omega^2&  \omega^2 & \omega   & 1   & \omega \\
1   &\omega^2& \omega & \omega^2    & \omega   &  1  
\end{array}\right)$$
is orthogonal\footnote{Unfortunately, we do not know of any slick way to verify this other than by brute force computation; this above example was discovered by trial and error.  It may have some algebraic interpretation, perhaps relating to a group of order 6 such as $\Z_6$ or $S_3$, but we were unable to discover such an interpretation.}, with all its entries equal to a cube root of unity.  By repeating
the above argument, we thus see that the six-element set
$\Omega'_0 := \{ e_1, \ldots, e_6 \}$ in $\Z_3^6$ has a spectrum $\Lambda'_0$, also with six elements, but does not tile $\Z_3^6$ since $6$ does not go evenly into $3^6$.

In fact one can descend from this example to $\Z_3^5$ by the following simple observation.  We may translate $\Omega'_0$ by $-e_1$ to create a new set
$$ \{ 0, e_2 - e_1, \ldots, e_6 - e_1 \};$$
this of course does not affect the property that $\Lambda'_0$ is a spectrum.  
This set is contained in the 5-dimensional space
$$ \Gamma := \{ (x_1, \ldots, x_6) \in \Z_3^6: x_1 + \ldots + x_6 = 0 \};$$
by projecting out the orthogonal complement $\Gamma^\perp := \{ (x,x,\ldots,x): x \in \Z_3 \}$ we may thus assume that $\Lambda'_0$ is also contained in $\Gamma$.  But $\Gamma$ is clearly isomorphic to $\Z_3^5$.  We have thus proved

\begin{theorem}\label{thm:finite-3}  There exists a set $\Omega''_0 \subset \Z_3^{5}$ of six elements, such that $l^2(\Omega''_0)$ admits an orthonormal basis of exponentials $\{ \frac{1}{(\# \Omega''_0)^{1/2}} e^{2\pi i (\xi_j \cdot x) / 3}: \xi_j \in \Lambda''_0 \}$ for some $\Lambda''_0 \subset \Z_3^5$, but such that $\Omega''_0$ does not tile $\Z_3^5$ by translations.
\end{theorem}

Note that the non-tiling property follows since $6$ does not go evenly 
into $3^5$.  A similar argument allows one to replace $\Z_2^{12}$ with $\Z_2^{11}$ in Theorem \ref{thm:finite}, but we will not need to do so here.  We also remark that the above theorem clearly also holds if $\Z_3^5$ is replaced by $\Z_3^n$ for any $n \geq 5$.

\section{The discrete model: failure of Fuglede in $\Z^5$}\label{sec:discrete}

We now modify the construction in the previous section to prove a discrete version of Theorem \ref{thm:main}, in the lattice $\Z^5$.  The dual group to this lattice is the torus $\R^5/\Z^5$; if $x \in \Z^5$ and $\xi \in \R^5/\Z^5$ we can define the expression $e^{2\pi i \xi \cdot x}$ in the obvious manner.

\begin{theorem}\label{thm:discrete}  There exists a non-empty finite subset $\Omega_1 \subset \Z^5$ such that $l^2(\Omega_1)$ admits an orthonormal basis of exponentials $\{ \frac{1}{(\# \Omega_1)^{1/2}} e^{2\pi i \xi_j \cdot x}: \xi_j \in \Lambda_1 \}$ for some $\Lambda_1 \subset \R^5 / \Z^5$, but such that $\Omega_1$ does not tile $\Z^5$ by translations.
\end{theorem}

\begin{proof}
Heuristically, this theorem follows immediately from Theorem \ref{thm:finite-3} by pulling $\Omega''_0$ back under the obvious homomorphism $\Z^5 \to \Z_3^5$.  Unfortunately this has the problem of making the pre-image of $\Omega''_0$ infinite; however this can be rectified by truncating this pre-image at a sufficiently large scale, and noting that boundary effects of the truncation will be negligible if the scale is large enough.

We turn to the details.
Let $\Omega''_0 \subset \Z_3^5$ and $\Lambda''_0 \subset \Z_3^5$ be the counterexample to Fuglede's conjecture in $\Z_3^5$ constructed in Theorem \ref{thm:finite-3}.  We need three large integers
$$ 1 \ll L \ll M \ll N,$$
for instance we may pick $L := 10^{10}$, $M := L^2$, and $N := M^2$.

We define $\Omega_1 \subset \Z^5$ to be the set
$$ \Omega_1 := \bigcup_{k \in [0,M)^5} (3k + \Omega''_0)$$
where we identify $\Z_3^5$ with the set $\{0,1,2\}^5 \subset \Z^5$ in the obvious (non-homomorphic!) manner, and $[0,M)^5$ is the discrete cube
$$ [0,M)^{5} := \{ (x_1, \ldots, x_5) \in \Z^5: 0 \leq x_j < M \hbox{ for all } j = 1, \ldots, 5 \}.$$
Observe that $\Omega_1$ is a finite set in $[0,3M)^{5}$ consisting of $6 M^{5}$ elements.  We define the spectrum $\Lambda_1 \subset \R^{5}/\Z^{5}$ in a similar fashion by
$$ \Lambda_1 := \bigcup_{l \in [0,M)^{5}} (\frac{l}{3M} + \frac{1}{3} \Lambda''_0),$$
where the homomorphism $\xi \mapsto \frac{1}{3}\xi$ from $\Z_3^{5}$ to $\{0+\Z,\frac{1}{3}+\Z, \frac{2}{3} + \Z\}^{5} \subset (\R/\Z)^{5}$ is defined in the obvious manner.  Note that $\Lambda_1$ is also a finite set consisting of $6 M^{5}$ elements.

We now verify that the set of exponentials $\{ \frac{1}{(\# \Omega_1)^{1/2}} e^{2\pi i \xi_j \cdot x}: \xi_j \in \Lambda_1 \}$ form an orthonormal basis of $l^2(\Omega_1)$.  The normalization property is obvious.  Since the number of exponentials equals the dimension of $l^2(\Omega_1)$, it will suffice to prove orthogonality, i.e. that
\begin{equation}\label{eq:ortho}
 \sum_{x \in \Omega''_0} e^{2\pi i (\xi - \xi') \cdot x} = 0 \hbox{ for all distinct } \xi,\xi' \in \Lambda_1.
\end{equation}
We write $\xi = \frac{l}{3M} + \frac{1}{3} \xi_0$ and $\xi' = \frac{l'}{3M} + \frac{1}{3} \xi'_0$ for some $l,l' \in [0,M)^{5}$ and $\xi_0, \xi'_0 \in \Lambda''_0$; since $\xi \neq \xi'$, we observe that at least one of $l \neq l'$ or $\xi_0 \neq \xi'_0$ must hold.  We similarly write $x = 3k + x_0$ where $k$ ranges over $[0,M)^{5}$ and $x_0$ ranges over $\Omega''_0$.  We can then rewrite the left-hand side of \eqref{eq:ortho} as
$$
\sum_{k \in [0,M)^{5}} \sum_{x_0 \in \Omega''_0} e^{2\pi i (\frac{l-l'}{3M} + \frac{1}{3} (\xi_0 - \xi'_0)) \cdot (3k + x_0)}.$$
We may expand this as
$$
\sum_{k \in [0,M)^{5}} \sum_{x_0 \in \Omega''_0} 
e^{2\pi i (l-l') \cdot k / M}
e^{2\pi i (l-l') \cdot x_0 / 3M}
e^{2\pi i (\xi_0-\xi'_0) \cdot x_0 / 3}
$$
since the dot product of $\frac{1}{3} (\xi_0 - \xi'_0)$ and $3k$ is an integer and hence negligible.  
 
The sum $\sum_{k \in [0,M)^{5}} e^{2\pi i (l-l') \cdot k / M}$ vanishes unless $l=l'$ (this is basically the Fourier inversion formula for $\Z/M\Z$).  Thus we may assume $l=l'$, and hence $\xi_0 \neq \xi'_0$.  But then the previous expression simplifies to
$$ M^{5} \sum_{x_0 \in \Omega''_0} e^{2\pi i (\xi_0-\xi'_0) \cdot x_0 / 3}$$
which vanishes since the frequencies in $\Lambda''_0$ were chosen to give an orthogonal basis of $l^2(\Omega''_0)$.  This proves the existence of a spectrum.

We now show that the set $\Omega_1$ does not tile $\Z^{5}$, if $L,M,N$ were chosen sufficiently large; this will be a volume packing argument that relies once again on the fact that $6$ does not go evenly into $3^{6}$.  Suppose for contradiction that we could find a subset $\Sigma \subset \Z^{5}$ such that the translates $\{ t + \Omega_1: t \in \Sigma \}$ tiled $\Z^{5}$.  

The idea is to exploit the intuitive observation that $\Omega_1$ has local density either equal to 0 or to $\frac{6}{3^{5}}$, except on the boundary.  To make this rigorous\footnote{Another approach, which is basically equivalent to this one, is to compute the convolution $\chi_{\Omega_1} * \frac{1}{L^n} \chi_{[0,L)^n}$ and observe that this is approximately equal to the function $\frac{6}{3^{5}} \chi_{[0,M)^n}$, which cannot tile $\Z^n$.}
we now use the other numbers $L, N$ chosen earlier. let $\Sigma_N := \Sigma \cap [0,N)^{5}$.  Observe that the sets $\{ t + \Omega_1: t \in \Sigma_N \}$ are disjoint and each have cardinality $6 M^{5}$, while their union contains the cube $[2M, N-2M)^{5}$, and is contained in the cube $[-2M, N+2M)^{5}$.  Since both of these cubes have a cardinality of $N^{5} + O(M N^4)$, we thus have the cardinality estimate
$$ 6 M^{5} \# \Sigma_N = N^{5} + O(M N^4).$$

Let $A$ denote the annulus
$$ A := [-5L, 2M+5L)^{5} - [5L, 2M-5L)^{5};$$
this $A$ represents the boundary effects of our restriction of $\Omega_1$ to $[0,M)^{5}$.  We shall now work in a region of the tiling where $A$ can be ignored.  Observe that $A$ has cardinality $\# A = O( M^4 L )$.  In particular, we have
$$ \# \bigcup_{t \in \Sigma_M} (t + A) \leq (\# \Sigma_N) (\# A)
= O( \frac{N^{5}}{6 M^{5}} M^4 L ) = O(\frac{L}{M}) N^{5}.$$
Thus, if $M$ is chosen sufficiently large with respect to $L$, we see that
$$ \# \bigcup_{t \in \Sigma_M} (t + A) < \# [\frac{N}{3}, \frac{2N}{3})^{5},$$
and thus we may find a point $x_0 \in [\frac{N}{3}, \frac{2N}{3})^{5}$ which is not contained in $t+A$ for any $t \in \Sigma_M$.

Fix this point $x_0$, and consider the quantity
$$ f(t) := \# \bigl( (t + \Omega_1) \cap (x_0 + [0,L)^{5}) \bigr)$$
defined for all $t \in \Sigma$.  Since the sets $\{ t + \Omega_1: t \in \Sigma \}$ tile $\Z^{5}$, we must have
\begin{equation}\label{eq:l-sum}
 \sum_{t \in \Sigma} f(t) = \#( x_0 + [0,L)^{5} ) = L^{5}.
\end{equation}
On the other hand, we shall shortly show that for every $t \in \Sigma$, we have either
\begin{equation}\label{eq:l-eq}
f(t) = 0 \hbox{ or } f(t) = L^{5} (\frac{6}{3^{5}} + O(\frac{1}{L})).
\end{equation}
Since $6$ does not go evenly into $3^{5}$, the estimates \eqref{eq:l-sum}, \eqref{eq:l-eq} will cause the desired contradiction if $L$ is chosen sufficiently large.

It remains to prove \eqref{eq:l-eq}.  We may assume of course that $f(t) \neq 0$; so that $t + \Omega_1$ intersects $x_0 + [0,L)^{5}$, which implies that \begin{equation}\label{eq:x0-loc}
x_0 \in t + [-L, 2M + L)^{5}.
\end{equation}
Since $x_0 \in [\frac{N}{3}, \frac{2N}{3}]^{5}$, this implies that $t \in [0,N]^{5}$, since $N$ is much larger than $L$ or $M$.  In particular, $t \in \Sigma_N$.  By construction of $x_0$, we thus see that $x_0 \not \in t + A$.  Combining this \eqref{eq:x0-loc} we see that
$$
x_0 \in t + [5L, 2M - 5L)^{5}.$$
In particular, we have
$$ x_0 + [0,L)^{5} \subset t + [4L, 2M-4L)^{5}.$$
If one then covers $x_0 + [0,L)^{5}$ by $\frac{1}{3^{5}} L^{5} + O(L^4)$ cubes of the form $2k + \{0,1\}^{5}$ for various $k \in \Z^{5}$ - and observe that all but $O(L^4)$ of these cubes will be strictly contained inside $x_0 + [0,L)^{5}$ - then each of these cubes will lie inside $t + [3L, 2M-3L)$, and thus will intersect $t+\Omega_1$ in exactly $\# \Omega''_0 = 5$ points.  Thus we have
$$ \# \bigl( (t+\Omega_1) \cap (x_0 + [0,L)^{5}) \bigr) = [\frac{1}{3^{5}} L^{5} + O(L^4)]5 + O(L^4)$$
which is \eqref{eq:l-eq} as desired.  This concludes the proof that $\Omega_1$ does not tile $\Z^{5}$, and we are done.
\end{proof}

\section{The continuous model: failure of Fuglede in $\R^{5}$}\label{sec:final}

We can now prove Theorem \ref{thm:main} by a standard transference argument.  Let $\Omega_1$ and $\Lambda_1$ be defined as in the previous section.  Then we simply set
$$ \Omega_2 := \Omega_1 + [0,1)^{5}$$
and
$$ \Lambda_2 := \Lambda_1 + \Z^{5}.$$
Since $\Lambda_1$ was a spectrum for $\Omega_1$, and $\Z^{5}$ is a spectrum for $[0,1)^{5}$, it is easy to see that $\Lambda_2$ is a spectrum for $\Omega_2$.  The proof that $\Omega_2$ does not tile $\R^{5}$ proceeds almost exactly the same as in the previous section, with the obvious changes that cubes such as $[0,N)^{5}$ should now be solid cubes in $\R^{5}$ rather than discrete cubes in $\Z^{5}$, and that cardinality $\# (A)$ of sets should mostly be replaced by Lebesgue measure $|A|$ instead (except when dealing with $\Sigma_N$, which remains discrete).  We omit the details.

The above argument can be extended from $\R^5$ to $\R^n$ for any $n \geq 5$, mainly becase Theorem \ref{thm:finite-3} also extends in this manner.

\end{document}